\documentstyle[amstex,amscd]{article}
\title{Mukai flops and derived categories}
\author{Yoshinori Namikawa}
\date{ }
\begin{document}
\maketitle

\begin{center}
{\bf Introduction}
\end{center}

Derived categories 
possibly give a new significant invariant for 
algebraic varieties. In particular, when 
the caninical line bundle is trivial, there 
are varieties which are not birationally 
equivalent, but have equivalent derived categories. 
The most typical example can be found in the 
original paper [Mu]. On the other hand, for such 
varieties, it is hoped that the birationally 
equivalence should imply the equivalence of 
derived categories. One of the important classes 
for testing this question is the class of complex 
symplectic manifolds (for varieties in other classes, 
see [B2, Ka 1, Ch]). A Mukai flop is a typical 
birational map between complex symplectic manifolds.      
  In this note, we shall prove that 
two smooth projective varieties 
of $\dim 2n$ connected by a Mukai 
flop have equivalent bounded derived 
categories of coherent sheaves.  
More precisely, let  
$X$ and $X^+$ be smooth projective 
varieties of 
dimension $2n$ such that there is a birational map 
$\phi: X - - \to X^+$ obtained as the Mukai 
flop along a subvariety $Y \subset X$ which 
is isomorphic to ${\bold P}^n$. 
By definition, $\phi$ is decomposed into the 
blowing-up $p: {\tilde X} \to X$ along $Y$ 
and the blowing down $p^+: \tilde X \to X^+$ 
contracting the 
$p$-exceptional divisor 
${\tilde Y}$ to the subvariety 
$Y^+\cong {\bold P}^n$ of $X^+$. 
On the other hand, there are 
birational morphisms 
$X \to {\bar X}$ and $X^+ \to 
{\bar X}$ which contract $Y$ and 
$Y^+$ to points respectively.  
We put ${\hat X}:= X\times_{{\bar X}}X^+$ 
and let $q: {\hat X} \to X$ and $q^+: {\hat X} 
\to X^+$ be the natural projections. 
Note that ${\hat X}$ is a normal crossing 
variety with two irreducible components 
${\tilde X}$ and ${\bold P}^n \times 
{\bold P}^n$.  
 
Let $D(X)$ (resp. $D(X^+)$) be the bounded 
derived category of coherent sheaves on 
$X$ (resp. $X^+$). 

We shall consider natural two functors 
$$\Phi := 
{\bold R}(p^+)_*\circ{\bold L}p^* : 
D(X) \to D(X^+)$$

$$\Psi := {\bold R}(q^+)_*\circ{\bold L}q^* : 
D(X) \to D(X^+).$$

In this note, we first take, as 
$X$ and $Y$, the  
${\bold P}^n$ bundle ${\bold P}(
{\cal O}_{{\bold P}^n} \oplus 
\Theta_{{\bold P}^n})$ over 
${\bold P}^n$ and its negative section. 
In this case,  
$\Phi$ is not a fully-faithful functor 
even when $n = 2$ (\S 2), but  
$\Psi$ is an equivalence of triangulated 
categories for any $n \geq 2$ (\S 3). 
In \S 4 we prove in a general 
case that $\Psi$ is an 
equivalence of triangulated categories 
(Theorem (4.4))
\footnote{Szendroi pointed out 
to me that a similar phenomenon to ours 
actually occurs in Calabi-Yau threefolds 
with type III contractions [Sz, Theorem 6.5].}.  
Theorem (4.4) holds for a Mukai flop in a 
more general sense. This generalization 
is done in \S 5. 
  
Our result should be compared with elementary 
flops studied in Bondal-Orlov [B-O, Theorem 
3.6]. For these 
flops,  
$\Phi$ becomes an equivalence 
of triangulated categories. But, 
in our case, since $\Phi \ne 
\Psi$, the proof there can not 
be directly applied. Instead, we 
shall use Bridgeland's criteria 
[B 1, Theorem 2.3], [B-K-K, Theorem 
2.4] for
exact functors of triangulated 
categories to be fully-faithful 
and to be an equivalence.  
A Mukai flop is of particular 
importance when $X$ is a  
complex irreducible symplectic manifold
; in this case, it is conjectured 
that the birationally equivalence 
implies the equivalence of 
bounded derived categories 
of coherent sheaves. 
This conjecture is also 
related with a question of 
Torelli type [Na, Question].
By a recent result of Wierzba and 
Wisniewski [W-W, Theorem 1.2](see 
also [W]), one 
can prove that, birationally equivalent 
complex projective symplectic 4-folds have 
equivalent derived categories (Corollary 
(4.5)). 

This is a revised version of 
the paper with the same title appeared 
in math.AG/0203287.  
Recently, Kawamata has independently proved 
the similar results in [Ka 2, $\S 5$]. 
In [Ka 2] the equivalence for Mukai 
flops is redued to the one for the 
elementary flops in [B-O]. 
But our proof here depend on direct calculations 
using an explicit model.       
\vspace{0.12cm}  

\begin{center}
{\bf \S 1.}
\end{center}
 
(1.1)  Fix an integer $n \geq 2$, and 
denote by $X$  
the projective space bundle 
${\bold P}({\cal O}_{{\bold P}^n} 
\oplus \Theta_{{\bold P}^n})$ over ${\bold P}^n$ 
in the sense of Grothendieck (cf. 
[Ha, p.162]). We put $M = {\bold P}^n$ and 
let $\pi: X \to M$ be the bundle map. 
Denote by ${\cal O}_X(1)$ 
the tautological line 
bundles on $X$.    
Let $Y \subset X$ be the section of $\pi$ defined by 
the surjection ${\cal O}_{{\bold P}^n} 
\oplus \Theta_{{\bold P}^n} \to 
{\cal O}_{{\bold P}^n}$.   

The normal bundle $N_{Y/X}$ is isomorphic to $\Omega^1_Y$. 
Let $p: {\tilde X} \to X$ be the blowing-up of $X$ along $Y$. 
The exceptional divisor ${\tilde Y}$ of $p$ is isomorphic 
to ${\bold P}(\Theta_{{\bold P}^n})$. ${\tilde Y}$ 
has a ${\bold P}^{n-1}$-bundle structure which 
is different from the bundle structure ${\tilde Y} \to Y$.  
One can contract ${\tilde Y} \subset {\tilde X}$ 
along this another ruling to $Y^+ \subset X^+$. 
We call this birational contraction map $p^+$. 
Let $\phi : X - - \to X^+$ be the birational 
map defined as the composite $p^+ \circ (p^{-1})$. 
$Y^+$ is isomorphic to ${\bold P}^n$ and $N_{Y^+/X^+} 
\cong \Omega^1_{Y^+}$.  \vspace{0.12cm}

{\bf Lemma (1.2)}. {\em $X^+$ is isomorphic 
to the projective space bundle 
${\bold P}({\cal O}_{{\bold P}^n} 
\oplus \Theta_{{\bold P}^n})$ over ${\bold P}^n$.}
\vspace{0.12cm}

{\em Proof}. Let $H$ be a hyperplane of $M$. Then 
$\pi^{-1}(H) = {\bold P}_H({\cal O}_H 
\oplus \Theta_M\vert_H)$. Let $V_H$ 
be the projective space subbundle of 
$\pi^{-1}(H) \to H$ defined by 
the surjection 
${\cal O}_H \oplus \Theta_M\vert_H 
\to {\cal O}_H \oplus N_{H/M}$. 
$V_H \to H$ has a negative section 
and the variety obtained by contracting 
it to a point becomes ${\bold P}^n$. 
Let $M^+$ be the dual projective 
space of $M$. Now $\{V_H\}_{H \in M^+}$ 
is a flat family of subvarieties of $X$. 
Let $V^+_H$ be the proper transform 
of $V_H$ by $\phi$. Then $\{V^+_H\}_{H 
\in M^+}$ is a flat family of subvarieties 
of $X^+$. Under the birational transformation 
$\phi$, the negative section of $V_H \to H$ 
is contracted to a point; hence each $V^+_H$ 
is isomorphic to ${\bold P}^n$. 
Moreover, one can check that $V^+_H$ 
and $V^+_{H'}$ are disjoint if $H \neq 
H'$. This implies that $X^+$ has 
a ${\bold P}^n$-bundle structure with 
the section $Y^+$. 
Let $X_{\infty} \subset X$ be the 
projective space subbundle defined 
by the surjection ${\cal O}_M 
\oplus \Theta_M \to \Theta_M$. 
$X_{\infty}$ is disjoint from $Y$. 
We see that the proper transform 
$X^+_{\infty}$ also becomes 
a projective space subbundle of 
$X^+$ disjoint from $Y^+$.  
Since $N_{Y^+/X^+} 
\cong \Omega^1_{Y^+}$, we 
have the result. 
\vspace{0.2cm} 

We denote by $\pi^+ : X^+ \to M^+$ the 
bundle structure introduced in (1.2) 
and denote by ${\cal O}_{X^+}(1)$ the 
tautological line bundle. \vspace{0.2cm}

{\bf Lemma (1.3)}. {\em The birational map 
$\phi$ induces an isomorpshism 
$$ \phi^* : {\mathrm Pic}(X^+) \to 
{\mathrm Pic}(X) $$ 
such that ${\phi}^*({\cal O}_{X^+}(1)) 
= {\cal O}_X(1)$ and ${\phi}^*(
(\pi^+)^*{\cal O}_{M^+}(1)) = 
{\cal O}_X(1)\otimes 
\pi^*{\cal O}_M(-1)$.} 
\vspace{0.15cm}

{\em Proof}. Since $\phi$ is an 
isomorphism in codimension 1, there 
is a natural isomorphism $\phi^*$ 
between ${\mathrm Pic}(X^+)$ and 
${\mathrm Pic}(X)$. By (1.2) we see 
that ${\phi}^*(
(\pi^+)^*{\cal O}_{M^+}(1)) = 
{\cal O}_X(1)\otimes 
\pi^*{\cal O}_M(a)$ for 
some $a \in {\bold Z}$. 
Since $h^0({\cal O}_X(1)\otimes 
\pi^*{\cal O}_M(a)) = 
h^0((\pi^+)^*{\cal O}_{M^+}(1)) 
= n + 1$, we conclude that $a = -1$.  

The birational transformation 
$\phi$ is symmetric with respect 
to $X$ and $X^+$; hence, we should 
have $(\phi^{-1})^*((\pi)^*{\cal O}_{M}(1)) 
= {\cal O}_{X^+}(1)\otimes 
(\pi^+)^*{\cal O}_{M^+}(-1)$ for 
${\phi^{-1}}^* : {\mathrm Pic}(X) 
\to {\mathrm Pic}(X^+)$.  

Then we have $\phi^* \circ (\phi^{-1})^* 
((\pi)^*{\cal O}_{M}(1)) = \phi^*({\cal 
O}_{X^+}(1)) \otimes {\cal O}_X(-1) \otimes 
\pi^*{\cal O}_M(1)$.  Since 
$\phi^* \circ (\phi^{-1})^* = id$,      
$\phi^*({\cal 
O}_{X^+}(1)) = {\cal O}_X(1)$.  
\vspace{0.2cm}

(1.4) For $k \in {\bold Z}$, let $D(X)_k$ be 
the full subcategory of $D(X)$ whose objects are 
of the form ${\cal O}_X(k)\otimes {\bold L}\pi^*(F)$ 
with $F \in D(M)$.  
As a triangulated category, $D(X)$ is generated by 
$D(X)_{-n}$, ..., $D(X)_{-1}$, $D(X)_0$ by 
[O 1, Theorem 2.6]. Since $D(M)$ is generated by 
the objects ${\cal O}_M(-n+1)$, ..., ${\cal O}_M$, 
and ${\cal O}_M(1)$, we see that $D(X)$ is generated 
by the set of objects $\Omega := \{ {\cal O}_X(j)\otimes 
\pi^*{\cal O}_M(k) \}$ with $-n \leq j \leq 0$ 
and $-n+1 \leq k \leq 1$. In particular, 
$\Omega$ is a spanning class for $D(X)$.  
Here, by definition of the spanning class, 
if an object $a \in D(X)$ satisfies ${\mathrm Hom}^i(\omega, 
a) = 0$ for all $\omega \in \Omega$ and for all $i \in {\bold Z}$, 
then $a \cong 0$. Similarly, if ${\mathrm Hom}^i(a, \omega) = 0$ 
for all $\omega \in \Omega$ and for all $i \in {\bold Z}$, 
then $a \cong 0$. 
\vspace{0.2cm} 

(1.5) Define functors $\Phi$ and $\Phi'$ as 

$$  \Phi := {\bold R}(p^+)_*{\bold L}p^*(\cdot): 
D(X) \to D(X^+). $$ 

$$  \Phi' := {\bold R}p_*({\bold L}(p^+)^*(\cdot)
\otimes{\omega_{{\tilde X}/X}}): 
D(X^+) \to D(X), $$ 
where $\omega_{{\tilde X}/X} := 
\omega_{\tilde X}\otimes p^*\omega^{-1}_X$. 
\vspace{0.15cm}

{\bf Lemma (1.6)}. {\em Assume 
that $j$ is an integer such that 
$-n \leq j \leq 0$.}  

(1) {\em When 
$-n+1 \leq k \leq 0$, 
$\Phi({\cal O}_X(j)\otimes 
\pi^*{\cal O}_M(k)) = 
{\cal O}_{X^+}(j+k)\otimes 
\pi^*{\cal O}_{M^+}(-k)$.}

{\em When $k = 1$, $\Phi({\cal O}_X(j)\otimes 
\pi^*{\cal O}_M(1)) = 
{\cal O}_{X^+}(j+1)\otimes 
\pi^*{\cal O}_{M^+}(-1)\otimes I_{Y^+}$,   
where $I_{Y^+}$ is the ideal sheaf of 
$Y^+$ in $X^+$.} 
\vspace{0.15cm}

(2) {\em When $-n+1 \leq k \leq 0$, 
$\Phi'\circ\Phi({\cal O}_X(j)\otimes 
\pi^*{\cal O}_M(k)) = {\cal O}_X(j)\otimes 
\pi^*{\cal O}_M(k)$.}   
\vspace{0.15cm}
   
{\em Proof}. (1): At first, we have 
${\cal O}_X(1)\vert_Y \cong 
{\cal O}_Y$, hence, for 
$y^+ \in Y^+$,  
$(p^*({\cal O}_X(j)\otimes   
\pi^*{\cal O}_M(k)))\vert_{{p^+}^{-1}(y^+)}$ is 
isomorphic to ${\cal O}_{{\bold P}^{n-1}}(k)$. 
Therefore, $p^*({\cal O}_X(j)\otimes 
\pi^*{\cal O}_M(k))$ is 
${p^+}_*$-acyclic for 
$k \geq -n+1$. If 
$k \leq 0$, then ${p^+}_*(p^*({\cal O}_X(j)\otimes 
\pi^*{\cal O}_M(k)))$ is a line bundle 
on $X^+$. If $k = 1$, then 
${p^+}_*(p^*({\cal O}_X(j)\otimes 
\pi^*{\cal O}_M(1)))$ is a line 
bundle outside $Y^+$ but not a 
line bundle on $Y^+$. 
  Finally, by (1.3) we conclude 
that, if $k \leq 0$, then      
$ \Phi({\cal O}_X(j)\otimes 
\pi^*{\cal O}_M(k)) = 
{\cal O}_{X^+}(j+k)\otimes 
\pi^*{\cal O}_{M^+}(-k)$, and that 
if $k = 1$, then $\Phi({\cal O}_X(j)\otimes 
\pi^*{\cal O}_M(1)) = 
{\cal O}_{X^+}(j+1)\otimes 
\pi^*{\cal O}_{M^+}(-1)\otimes I_{Y^+},$  
where $I_{Y^+}$ is the ideal sheaf of 
$Y^+$ in $X^+$.  \vspace{0.15cm}

(2): For $y \in Y$, $(p^+)^*({\cal O}_{X^+}
(j+k)\otimes(\pi^+)^*{\cal O}_{M^+}(-k))
\otimes \omega_{{\tilde X}/X}\vert_{p^{-1}(y)}$ 
is isomorphic to 
${\cal O}_{{\bold P}^{n-1}}(-k-n+1)$ because 
$\omega_{{\tilde X}/X}\vert_{p^{-1}(y)} 
\cong {\cal O}_{{\bold P}^{n-1}}(-n+1)$. 
By the same argument as (1), we conclude that 
$$ \Phi'({\cal O}_{X^+}
(j+k)\otimes(\pi^+)^*{\cal O}_{M^+}(-k)) 
= {\cal O}_X(j)\otimes 
\pi^*{\cal O}_M(k).$$

\vspace{0.2cm}

(1.7) We shall construct a locally 
free resolution 
of the ideal sheaf $I_{Y^+}$. 

Put $E := {\cal O}_{M^+}\oplus\Theta_{M^+}$. 
Let 
$$ 0 \to {\cal O}_{X^+}(-1) \to 
(\pi^+)^*E^* 
\to Q \to 0$$ 
be the universal exact sequence, where 
$E^*$ is the dual of $E$. 
Let us consider the fiber product 
$X^+ \times_{M^+}X^+$ and denote by 
$p_1: X^+ \times_{M^+}X^+ \to X^+$ and 
$p_2: X^+ \times_{M^+}X^+ \to X^+$ the 
first projection and the second projection 
respectively. One has an isomorphism ([O 1, p.136]) 
$$ H^0(X^+ \times_{M^+}X^+, {p_1}^*{\cal O}_{X^+}(1)
\otimes{p_2}^*Q) \cong H^0(M^+, E\otimes E^*).$$    
Corresponding to the identity map $id: E \to E$, 
one has a section $s \in  
H^0(X^+ \times_{M^+}X^+, {p_1}^*{\cal O}_{X^+}(1)
\otimes{p_2}^*Q)$. Then the set of zeros of $s$ 
coincides with the diagonal $\Delta \subset 
X^+ \times_{M^+}X^+$ (cf. [O 1, p.136]). Hence,   
we have the following Koszul complex which 
becomes a resolution of 
the ${\cal O}_{X^+ \times_{M^+}X^+}$ module  
${\cal O}_{\Delta}$: \vspace{0.12cm}

$ 0 \to {\wedge}^n({p_1}^*{\cal O}_{X^{+}}(-1)\otimes 
{p_2}^*Q^*) \to 
{\wedge}^{n-1}({p_1}^*{\cal O}_{X^{+}}(-1)\otimes 
{p_2}^*Q^*) \to ... \to $

${p_1}^*{\cal O}_{X^{+}}(-1)\otimes 
{p_2}^*Q^* \to {\cal O}_{{X^+}\times_{M^+}X^+} 
\to {\cal O}_{\Delta} \to 0.$ 

\vspace{0.12cm}

By the projection $p_2 : X^+ \times_{M^+}X^+ \to X^+$, 
each term becomes a flat ${\cal O}_{X^+}$ module. 
Thus, by taking the tensor product of the complex 
with the ${\cal O}_{X^+}$ module ${\cal O}_{Y^+}$, 
we get the exact sequence \vspace{0.12cm}

$(*)$

$ 0 \to {\wedge}^n({\cal O}_{X^+}(-1)\otimes 
(\pi^+)^*\Theta_{M^+}) \to 
{\wedge}^{n-1}({\cal O}_{X^+}(-1)\otimes 
(\pi^+)^*\Theta_{M^+}) \to ... \to $

${\cal O}_{X^+}(-1)\otimes 
(\pi^+)^*\Theta_{M^+})   
\to I_{Y^+} \to 0.$

\vspace{0.12cm}

Here we have used the fact that 
$\mathrm{Ker}[{\cal O}_{X^+} \to 
{\cal O}_{Y^+}] = I_{Y^+}$ and 
$Q^*\vert_{Y^+} \cong \Theta_{Y^+}$.     
\vspace{0.15cm}

\begin{center}
{\bf \S 2.}
\end{center}

In this section, we observe that 
the functor $\Phi$ is not fully 
faithful. \vspace{0.12cm}

{\bf Lemma (2.1)}. {\em Assume that 
$n = 2$, that is, $\dim X = \dim X^+ 
= 4$. Then 
${\mathrm Ext}^2(\Phi(\pi^*{\cal O}_M(1)), 
\Phi(\pi^*{\cal O}_M(1)) \neq 0$.}  
\vspace{0.15cm}

{\em Proof}. By (1.6) we have \vspace{0.12cm}

${\mathrm Ext}^2(\Phi(\pi^*{\cal O}_M(1)), 
\Phi(\pi^*{\cal O}_M(1)) =$ 
 
${\mathrm Ext}^2({\cal O}_{X^+}(1)\otimes 
\pi^*{\cal O}_{M^+}(-1)\otimes I_{Y^+}, 
{\cal O}_{X^+}(1)\otimes 
\pi^*{\cal O}_{M^+}(-1)\otimes I_{Y^+}) =$ 

${\mathrm Ext}^2(I_{Y^+}, I_{Y^+}).$  
\vspace{0.12cm}

By (*) in \S 1, we have an exact sequence 
$$  0 \to {\wedge}^2({\cal O}_{X^+}(-1)\otimes 
(\pi^+)^*\Theta_{M^+}) \to 
{\cal O}_{X^+}(-1)\otimes 
(\pi^+)^*\Theta_{M^+}   
\to I_{Y^+} \to 0. $$

This exact sequence yields the exact 
sequence \vspace{0.12cm}

$ {\mathrm Ext}^1({\cal O}_{X^+}(-1)\otimes 
(\pi^+)^*\Theta_{M^+}, I_{Y^+}) \to 
{\mathrm Ext}^1({\wedge}^2({\cal O}_{X^+}(-1)\otimes 
(\pi^+)^*\Theta_{M^+}, I_{Y^+})$ 

$\to {\mathrm Ext}^2(I_{Y^+}, I_{Y^+}) \to 
{\mathrm Ext}^2({\cal O}_{X^+}(-1)\otimes 
(\pi^+)^*\Theta_{M^+}, I_{Y^+}).$  
\vspace{0.12cm}

Note that, since ${\cal O}_{X^+}(-1)\otimes 
(\pi^+)^*\Theta_{M^+}$ and 
${\wedge}^2({\cal O}_{X^+}(-1)\otimes 
(\pi^+)^*\Theta_{M^+}$ are both locally free 
sheaves, the first, second and fourth terms 
are identified with suitable cohomology 
groups. 
Then, by a direct calculation one can check that 
the first and fourth term vanish but the second 
one is not zero. \vspace{0.15cm}
  
{\bf Corollary (2.2)}. {\em When $n = 2$, 
the functor $\Phi : D(X) \to D(X^+)$ is 
not fully faithful.} \vspace{0.15cm}

{\em Proof}. One has \vspace{0.12cm}

${\mathrm Hom}(\pi^*{\cal O}_M(1), 
\pi^*{\cal O}_M(1)[2]) = {\mathrm Ext}^2
(\pi^*{\cal O}_M(1), 
\pi^*{\cal O}_M(1)) =$ 

$H^2(X, {\cal O}_X) = 0$.   
\vspace{0.15cm}

On the other hand, by (2.1) 
\vspace{0.12cm}

${\mathrm Hom}(\Phi(\pi^*{\cal O}_M(1)), 
\Phi(\pi^*{\cal O}_M(1)[2])) = 
{\mathrm Hom}(\Phi(\pi^*{\cal O}_M(1)), 
\Phi(\pi^*{\cal O}_M(1))[2]) =$ 

${\mathrm Ext}^2
(\Phi(\pi^*{\cal O}_M(1)), 
\Phi(\pi^*{\cal O}_M(1))) \neq 0.$  
\vspace{0.2cm}

We shall give 
a more intrinsic proof to 
Corollary (2.2). 
\vspace{0.15cm}

For this we need a lemma. 
\vspace{0.12cm}

{\bf Lemma (2.3)}. {\em Notation being the same as 
$\S 1$. Then, ${\mathrm Ext}^i({\cal O}_{Y^+}, 
{\cal O}_{Y^+}) = {\bold C}$ if $i$ is an even integer 
with $0 \leq i \leq 2n$, and otherwise 
${\mathrm Ext}^i({\cal O}_{Y^+}, 
{\cal O}_{Y^+}) = 0.$} 
\vspace{0.15cm} 

{\em Proof}. By the exact sequence 
$$ 0 \to I_{Y^+} \to {\cal O}_{X^+} \to 
{\cal O}_{Y^+} \to 0 $$
one has the exact sequence
\vspace{0.15cm}
 
$ \underline{{\mathrm Hom}}({\cal O}_{Y^+}, 
{\cal O}_{Y^+}) \to 
\underline{{\mathrm Hom}}({\cal O}_{X^+}, 
{\cal O}_{Y^+}) \to$ 

$\underline{{\mathrm Hom}}(I_{Y^+}, 
{\cal O}_{Y^+}) \to \underline{{\mathrm Ext}}^1({\cal O}_{Y^+}, 
{\cal O}_{Y^+}) \to 0.$
\vspace{0.15cm}

Since the first map is an isomorphism, 
$\underline{{\mathrm Hom}}(I_{Y^+}, 
{\cal O}_{Y^+}) \cong \underline{{\mathrm Ext}}^1({\cal O}_{Y^+}, 
{\cal O}_{Y^+})$. Here the left hand side is the normal 
bundle $N_{Y^+/X^+}$ which is isomorphic to $\Omega^1_{Y^+}$. 
Hence, we have 
$$ \underline{{\mathrm Ext}}^1({\cal O}_{Y^+}, 
{\cal O}_{Y^+}) \cong \Omega^1_{Y^+}. $$ 

Now, by taking $\underline{{\mathrm Hom}}(\cdot , 
{\cal O}_{Y^+})$ of the Koszul complex (*) in (1.7) we 
see that 
$$ \underline{{\mathrm Ext}}^i({\cal O}_{Y^+}, 
{\cal O}_{Y^+}) \cong 
\wedge^i \underline{{\mathrm Ext}}^1({\cal O}_{Y^+}, 
{\cal O}_{Y^+}) \cong \Omega^i_{Y^+}. $$ 

Therefore, in the spectral sequence 
$$ E^{p,q}_2 := H^p(X^+, \underline{{\mathrm Ext}}^q
({\cal O}_{Y^+}, {\cal O}_{Y^+})) => {\mathrm Ext}^{p+q}
({\cal O}_{Y^+}, {\cal O}_{Y^+}),$$
$E^{p,q}_2 = 0$ if $p \neq q$. 
This implies that this spectral sequence degenerates 
at $E_2$-terms; hence we have 
the result.  
\vspace{0.15cm} 

(2.4)({\em Another proof of (2.2)}): 
 For simplicity, 
put $L := \pi^*{\cal O}_M(1)$ 
and $L^+ := {\cal O}_{X^+}
(1)\otimes (\pi^+)^*{\cal O}_{M^+}
(-1)$. Then, $\Phi (L) = L^+ 
\otimes I_{Y^+}$. 
We have   
${\mathrm Ext}^2(L, L) 
= H^2(X, {\cal O}_X)$. 
On the other hand, 
${\mathrm Ext}^2(L^+ 
\otimes I_{Y^+}, 
L^+ \otimes I_{Y^+}) 
= {\mathrm Ext}^2(I_{Y^+}, 
I_{Y^+})$. 
Take ${\mathrm Hom}(., 
I_{Y^+})$ of the exact sequence 
$ 0 \to I_{Y^+} 
\to {\cal O}_{X^+} 
\to {\cal O}_{Y^+} \to 0$. 
Here  
${\mathrm Ext}^2({\cal O}_{Y^+}, 
I_{Y^+}) = 0$ because, in 
the exact sequence 
$$ {\mathrm Ext}^1({\cal O}_{Y^+}, 
{\cal O}_{Y^+}) \to 
{\mathrm Ext}^2({\cal O}_{Y^+}, 
I_{Y^+}) \to {\mathrm Ext}^2
({\cal O}_{Y^+}, {\cal O}_{X^+}) =0,$$ 
the first term vanish by Lemma (2.3). 
Moreover, we calculate \vspace{0.12cm}
 
${\mathrm 
Ext}^3({\cal O}_{X^+}, I_{Y^+}) 
= H^3(X^+, I_{Y^+}) \cong 
H^3(X^+, {\cal O}_{X^+}) = 0.$ 
\vspace{0.12cm}
   
Therefore, we have an exact sequence 
$$ 0 \to 
{\mathrm Ext}^2({\cal O}_{X^+}, 
I_{Y^+}) \to 
{\mathrm Ext}^2(I_{Y^+}, I_{Y^+}) 
\to {\mathrm Ext}^3({\cal O}_{Y^+}, 
I_{Y^+}) \to 0. $$ 
The first term is isomorphic 
to $H^2(X^+, I_{Y^+}) 
\cong H^2(X^+, {\cal O}_{X^+})$. 
Note that, $h^2(X^+, 
{\cal O}_{X^+}) = 
h^2(X, {\cal O}_X)$.  
The third term is isomorphic to 
${\mathrm Ext}^2({\cal O}_{Y^+}, 
{\cal O}_{Y^+}).$ 
By Lemma (2.3) this is a
one dimensional ${\bold C}$-vector 
space. Therefore, 
$$ \dim {\mathrm Ext}^2(L, L) = 
\dim {\mathrm Ext}^2(\Phi (L), 
\Phi (L)) - 1.$$

\vspace{0.2cm}

\begin{center}
{\bf \S 3.}
\end{center}

There are 
birational morphisms 
$X \to {\bar X}$ and $X^+ \to 
{\bar X}$ which contract $Y$ and 
$Y^+$ to points respectively.  
We put ${\hat X}:= X\times_{{\bar X}}X^+$ 
and let $q: {\hat X} \to X$ and $q^+: {\hat X} 
\to X^+$ be the natural projections. 
Note that ${\hat X}$ is a normal crossing 
variety with two irreducible components 
${\tilde X}$ and $Y \times 
Y^+$.  
Define a functor $\Psi$ as  
$$\Psi := {\bold R}(q^+)_*\circ{\bold L}q^* : 
D(X) \to D(X^+).$$

{\bf Theorem (3.1)}. {\em The functor 
$\Psi$ is an equivalence of triangulated 
categories.} 
\vspace{0.15cm}

The remainder consists of the proof of this 
theorem. 

(3.1){\em Outline of the 
Proof}:  We define a suitable 
spanning class $\Omega'$ for $D(X)$. 
We shall prove that, for any $a$, $b 
\in \Omega$, ${\mathrm Hom}^i(a,b) 
\cong {\mathrm Hom}^i(\Psi(a), \Psi(b))$ 
for all $i \in {\bold Z}$. Then, by 
[B 1, Theorem 2.3] we see that $\Psi$ is 
fully-faithful. 

Finally, apply Theorem  
2.4. in [B-K-R] to $\Psi$ 
and the spanning class 
$\Omega'$ to 
conclude that $\Psi$ is an equivalence. 
\vspace{0.12cm}

(3.2) For $k \in {\bold Z}$, let $D(X)_k$ be 
the full subcategory of $D(X)$ whose objects are 
of the form ${\cal O}_X(k)\otimes {\bold L}\pi^*(F)$ 
with $F \in D(M)$.  
As a triangulated category, $D(X)$ is generated by 
$D(X)_{-n}$, ..., $D(X)_{-1}$, $D(X)_0$ by 
[O, Theorem 2.6]. Since $D(M)$ is generated by 
the objects ${\cal O}_M(-n)$, ..., ${\cal O}_M$, 
 we see that $D(X)$ is generated 
by the set of objects $\Omega' := \{ {\cal O}_X(j)\otimes 
\pi^*{\cal O}_M(k) \}$ with $-n \leq j \leq 0$ 
and $-n \leq k \leq 0$. In particular, 
$\Omega'$ is a spanning class of $D(X)$. 

In the remainder, $j$ is an integer 
such that $-n \leq j \leq 0$, and 
$k$ is an integer such that $-n 
\leq k \leq 0$. 
\vspace{0.15cm}

(3.3) Let $f : Y \times Y^+ \to Y$ 
and $f^+: Y \times Y^+ \to Y^+$ be the 
natural projections, respectively. 
Note that ${\tilde Y}$ is a subvariety 
of $Y \times Y^+$. 
For an ${\cal O}_Y$ module $F$ and 
for an ${\cal O}_{Y^+}$ module $G$, write 
$F \odot G$ for $f^*F \otimes (f^+)^*G$. 
We write $F \odot_{\tilde Y}G$ for 
$f^*F \otimes (f^+)^*G\vert_{\tilde Y}$.   

We have an exact sequence (cf. [Fr]) 
$$ 0 \to {\cal O}_{\hat X} \to 
{\cal O}_{\tilde X}\oplus {\cal O}_{Y\times 
Y^+} \to {\cal O}_{\tilde Y} \to 0.$$   

Taking the tensor product of this sequence 
with $q^*({\cal O}_X(j)\otimes\pi^*{\cal O}_M(k))$, 
we get the exact sequence 

$ 0 \to q^*({\cal O}_X(j)\otimes\pi^*{\cal O}_M(k))
\to p^*({\cal O}_X(j)\otimes\pi^*{\cal O}_M(k))\oplus 
({\cal O}_Y(k)\odot {\cal O}_{Y^+})$ 

$\to {\cal O}_Y(k)\odot_{\tilde Y}{\cal O}_{Y^+} \to 0.$  

Here we have used the fact that ${\cal O}_X(1)\vert_Y 
\cong {\cal O}_Y$. \vspace{0.15cm}

When $k \geq -n+1$, the second (non-zero) term 
and the third one are both $q^+_*$-acyclic. 
Moreover, the maps 
$\alpha_k: q^+_*({\cal O}_Y(k)\odot {\cal O}_{Y^+})  
\to q^+_*({\cal O}_Y(k)\odot_{\tilde Y}{\cal O}_{Y^+})$ 
are surjective; hence the first (non-zero) term of the 
exact sequence is 
also $q^+_*$-acyclic. 
When $k = -n$, $R^{n-1}q^+_*
(p^*({\cal O}_X(j)\otimes\pi^*{\cal O}_M(-n)) 
\cong 
R^{n-1}q^+_*({\cal O}_Y(-n)\odot_{\tilde Y}{\cal O}_{Y^+})$. 
Moreover, ${\cal O}_Y(-n)\odot {\cal O}_{Y^+}$
is $q^+_*$-acyclic.  So, in this case, the first (non-zero) 
term of the exact sequence is $q^+_*$-acyclic. 
Therefore, when $-n \leq k \leq 0$, we have  

$\Psi({\cal O}_X(j)\otimes \pi^*{\cal O}_M(k)) 
= q^+_*q^*({\cal O}_X(j)\otimes \pi^*{\cal O}_M(k)) 
= p^+_*p^*({\cal O}_X(j)\otimes \pi^*{\cal O}_M(k)).$  
\vspace{0.12cm}
By the same calculation as the proof of Lemma (1.6), (1), 
$p^+_*p^*({\cal O}_X(j)\otimes \pi^*{\cal O}_M(k) 
= {\cal O}_{X^+}(j+k)\otimes (\pi^+)^*{\cal O}_{M^+}(-k)$ 
because $k \leq 0$.  As a consequence, we have  
$$ \Psi({\cal O}_X(j)\otimes \pi^*{\cal O}_M(k)) 
\cong {\cal O}_{X^+}(j+k)\otimes (\pi^+)^*{\cal O}_{M^+}(-k).$$ 
\vspace{0.15cm}

{\bf Lemma (3.4)} {\em Let $\pi: X \to M$ be the same as 
$\S 1$. Let $l$ and $m$ be integers such that 
$-n \leq l \leq n$ and $-n \leq m \leq n$. Then we have:} 
\vspace{0.15cm}

(1) $H^i(X, {\cal O}_X(l)\otimes {\cal O}_M(m)) = 0$ 
for $i > 0$.  \vspace{0.12cm}

(2) $H^i(X, {\cal O}_X(l+m)\otimes {\cal O}_M(-m)) = 0$ 
for $i > 0$. \vspace{0.15cm}

{\em Proof}. (1): Since $l \geq -n$, $H^i(X, 
{\cal O}_X(l)\otimes {\cal O}_M(m)) = H^i(M, 
\pi_*{\cal O}_X(l)\otimes {\cal O}_M(m))$.  
We only have to consider the case where $l \geq 0$.  
In this case $\pi_*{\cal O}_X(l)\otimes{\cal O}_M(m) 
\cong \mathrm{Sym}^l({\cal O}_M \oplus \Theta_M)\otimes 
{\cal O}_M(m)$. Now (1) follows from the fact that 
$m \geq -n$.  \vspace{0.15cm}

(2): When $l+m \geq -n$, $H^i(X, 
{\cal O}_X(l+m)\otimes {\cal O}_M(-m)) 
= H^i(M, \pi_*({\cal O}_X(l+m)\otimes {\cal O}_M(-m)))$. 
Since $-m \geq -n$, (2) is verified by the same 
argument as (1).  
When $l+m < -n$, we use the Serre duality.  
Note that $\omega_X \cong {\cal O}_X(-n-1)$. 
We have to prove that $H^{2n-i}(X, {\cal O}_X(-l-m 
-n-1)\otimes \pi^*{\cal O}_M(m)) = 0$ when  
$i > 0$.  
Since $-l-m-n-1 \geq 0$, we must show 
that $H^{2n-i}(M, \pi_*{\cal O}_X(-l-m 
-n-1)\otimes{\cal O}_M(m)) = 0$  
Since $m \geq -n$, this is true for 
$i \ne 2n$. When $i = 2n$, this is proved 
by a case-by-case checking. For example, when 
$l+m = -2n$, we must have $l = -n$ and 
$m = -n$. But, then $\pi_*{\cal O}_X(n-1) 
\otimes {\cal O}_M(-n) = \mathrm{Sym}^{n-1}
({\cal O}_M \oplus \Theta_M)\otimes {\cal O}_M(-n)$ 
has no global sections. \vspace{0.15cm}

{\bf Proposition (3.5)} {\em Let $a$ and $b$ be 
elements of $\Omega'$. Then, for all $i \in {\bold 
Z}$, ${\mathrm Hom}^i(a,b) \cong {\mathrm Hom}^i(\Psi(a), 
\Psi(b))$.} \vspace{0.15cm}

{\em Proof}. By Lemma (3.4),(1), we have 
${\mathrm Hom}^i(a,b) = 0$ for $i > 0$. By (3.3) 
and Lemma (3.4),(2), we have ${\mathrm Hom}^i(\Psi(a), 
\Psi(b)) = 0$ for $i > 0$. Let $a$ and $b$ be reprenented 
by the line bundles $L$ and $L'$ on $X$. Then 
$\Psi(a)$ and $\Psi(b)$ are represented by 
the proper transform $\phi(L)$ and $\phi(L')$ 
where $\phi : X --\to X^+$ is the Mukai flop. 
Therefore, ${\mathrm Hom}^0(a,b) \cong 
{\mathrm Hom}^0(\Psi(a), \Psi(b))$.     
\vspace{0.2cm}

(3.6).  By (3.5), for 
any $a$, $b \in \Omega'$, and 
for all $i \in {\bold Z}$, 
${\mathrm Hom}^i(a, b) 
\cong {\mathrm Hom}^i(\Psi (a), 
\Psi (b))$. 
$\Psi$ has a left adjoint and a right 
adjoint; in fact, 
let $p_1 : X\times X^+ \to X$ 
and $p_2 : X\times X^+ \to X^+$ be the 
first and second projections. Then 
$$\Psi (\cdot) = {\bold R}(p_2)_*({\cal O}
_{\hat X}\otimes^{\bold L}_
{{\cal O}_{X \times X^+}}p_2^*(\cdot)).$$ 
Now, put $Q:= {\bold R}{\mathrm Hom}_{D(X \times 
X^+)}({\cal O}_{\hat X}, {\cal O}_{X \times X^+})
\otimes (p_2)^*\omega_{X^+}[2n]$. 
Then, $$\Psi' (\cdot) := 
{\bold R}(p_1)_*(Q \otimes^{\bold L}
(p_2)^*(\cdot)): D(X^+) \to D(X)$$ 
becomes the left adjoint of $\Psi$ 
(cf. [B-K-R, \S 6, Step 1]).        
Let $S_X: D(X) \to D(X)$ be a Serre functor 
defined as 
$S_X(\cdot) = \cdot \otimes \omega_X [2n]$. 
Similary, we define $S_{X^+}: D(X^+) \to 
D(X^+)$.  Then 
$S_X \circ \Psi' \circ S_{X^+}^{-1}$ becomes 
a right adjoint of $\Psi$.  
Since $\Omega$ is a 
spanning class, we conclude that 
$\Psi$ is fully-faithful by [B 1, Theorem 2.3].  
By [B-K-R, Theorem 2.4], in order to prove that $\Psi$ is 
an equivalence of triangulated categories, 
we have to check that $\Psi \circ S_X (c) = 
S_{X^+} \circ \Psi (c)$ for all $c \in \Omega'$. 
Since $\omega_X \cong 
{\cal O}_X(-n-1)$, $\Psi \circ S_X (c) = 
\Psi (c) \otimes {\cal O}_{X^+}(-n-1)[2n]$.  
This coincides with $S_{X^+}\circ \Psi (c)$ 
because $\omega_{X^+} \cong 
{\cal O}_{X^+}(-n-1)$.  
\vspace{0.15cm}

\begin{center}
{\bf \S 4. } 
\end{center}

The result in the previous section was concerned 
with a special example. In this section we prove 
that the same holds in a general situation. 

(4.1)  Let $Z$ and $Z^+$ be birationally equivalent smooth 
projective varieties of dimension $2n$ with $n \geq 2$. 
Assume that there are subvarieties $W \subset Z$ 
and $W^+ \subset Z^+$ which are isomorphic to 
${\bold P}^n$. Assume that $N_{W/Z} \cong 
\Omega^1_{{\bold P}^n}$, $N_{W^+/Z^+} 
\cong \Omega^1_{{\bold P}^n}$ and that $Z^+$ 
(resp. $Z$) is  
the Mukai flop of $Z$ (resp. $Z^+$) 
along $W$ (resp. $W^+$).  
Let $Z \to \overline{Z}$ and 
$Z^+ \to \overline{Z}$ be the birational 
morphisms which contract $W$ and $W^+$ 
to points respectively. We put ${\hat Z} := 
Z \times_{\overline{Z}}Z^+$, and denote by 
$u: {\hat Z} \to Z$ and $u^+ : {\hat Z} \to Z^+$ 
the natural projections. We define the functor 
$$\Psi_Z : = {\bold R}(u^+)_*\circ {\bold L}u^* : 
D(Z) \to D(Z^+). $$ 

(4.2)  Let $X$, $X^+$ and others be the same as the 
previous sections. Define a set of objects 
$\Omega'' := \{{\cal O}_x\}_{x \in X}$, where 
${\cal O}_x$ is the structure sheaf of a point 
$x$.  Then $\Omega''$ becomes a spanning 
class for $D(X)$. Since the functor $\Psi : 
D(X) \to D(X^+)$ is an equivalence of triangulated 
categories, for any points $x$, $x' \in X$, 
${\mathrm Hom}^i({\cal O}_x, {\cal O}_{x'}) 
\cong {\mathrm Hom}^i(\Psi ({\cal O}_x), 
\Psi ({\cal O}_{x'}))
$, $i \in {\bold Z}$.   
Let $X_Y$ (resp. $X^+_{Y^+}$) be the 
formal completion of 
$X$ (resp. $X^+$) along $Y$ (resp. $Y^+$). 
Now assume that $x \in Y$ and $x' \in Y$. 
Since the supports of ${\cal O}_x$ 
and ${\cal O}_{x'}$ (resp. $\Psi ({\cal O}_x)$ 
and $\Psi ({\cal O}_{x'})$ ) are in $Y$ 
(resp. $Y^+$), we have 

$$  {\mathrm Hom}^i({\cal O}_x, {\cal O}_{x'}) 
= {\mathrm Hom}^i_{D(X_Y)}({\cal O}_x, {\cal O}_{x'}), $$ 
and 
$$  {\mathrm Hom}^i(\Psi ({\cal O}_x), 
\Psi ({\cal O}_{x'}))
= {\mathrm Hom}^i_{D(X^+_{Y^+})}(\Psi ({\cal O}_x), 
\Psi ({\cal O}_{x'})). $$ 

(4.3)  Let $Z_W$ (resp. $Z^+_{W^+}$) be 
the formal completion of $Z$ (resp. $Z^+$) 
along $W$ (resp. $W^+$). We have 
$Z_W \cong X_Y$  and 
$Z^+_{W^+} \cong X^+_{Y^+}$. 
By these identifications, we regard 
a point $x \in Y$ as a point on $W$. 
Now, by (4.2), if $x \in W$ and $x' \in W$, 
then we have 
$$ {\mathrm Hom}^i({\cal O}_x, {\cal O}_{x'}) 
\cong {\mathrm Hom}^i(\Psi_Z ({\cal O}_x), 
\Psi_Z ({\cal O}_{x'})). $$ 
When one of $x \in Z$ and $x' \in Z$ is not in 
$W$, it is clear that both spaces are isomorphic. 
Therefore, by [B 1, Theorem 2.3] $\Psi_Z$ is 
fully-faithful.  Let $\Omega''_Z$ be the spanning 
class for $D(Z)$ consisting of the objects 
${\cal O}_x$ with $x \in Z$, and let $S_Z$ (resp. 
$S_{Z^+}$) be the Serre functor defined in (3.10). 
Since $\omega_{Z^+}$ is trivial along 
$W^+$,  it is easily checked that  
$\Psi_Z \circ S_Z (c) = 
S_{Z^+} \circ \Psi_Z (c)$ for all $c \in \Omega''_Z$. 
Hence, by [B-K-R, Theorem 2.4], 
$\Psi_Z$ is an equivalence of triangulated categories.  
As a consequence, we have proved that 
\vspace{0.15cm}

{\bf Theorem (4.4)}. {\em Let $Z$, $Z^+$ and 
$\Psi_Z$ be the same as (4.1). Then 
$$ \Psi_Z : D(Z) \to D(Z^+) $$ 
is an equivalence of triangulated categories. }    
\vspace{0.12cm}

{\bf Corollary (4.5)}. {\em Let $Z$ and $Z'$ be 
birationally equivalent, complex 
projective symplectic 4-folds. Then 
$D(Z)$ and 
$D(Z')$ are equivalent.} 
\vspace{0.12cm}

{\em Proof}. By Wierzba and Wisniewski [W-W, 
Theorem 1.2](see also [W]), $Z$ and $Z'$ are 
connected by a finite sequence of 
Mukai flops. Theorem (4.4) together 
with this implies that $D(Z)$ and 
$D(Z')$ are equivalent.  
\vspace{0.15cm}

\begin{center}
{\bf \S 5.}
\end{center}
  
Theorem (4.4) holds for 
a Mukai flop in a more general sense. Namely, 
let $Z$ be a smooth projective 
variety of dim $2n$. Let $f: Z \to {\bar Z}$ be a 
projective birational morphism which contarcts a 
smooth subvariety $W \subset Z$ of dim $n+r$ 
($0 \leq r \leq n-2$) to a smooth subvariety $S 
\subset {\bar Z}$ of dim $2r$. Assume that $W \to 
S$ is a ${\bold P}^{n-r}$ bundle and assume that, for 
all $s \in S$, $N_{W/Z}\vert_{W_s} 
\cong \Omega^1_{W_s}$. Then we can perform Mukai 
flops in a family to get a new variety $Z^+$. Let 
us consider $Z \times_{{\bar Z}}Z^+$ and let $\Psi_Z: 
D(Z) \to D(Z^+)$ be the corresponding functor. Then 
we have the following. 
\vspace{0.15cm}

{\bf Theorem (5.1)}. {\em $\Psi_Z$ is an equivalence of triangulated 
categories.}  \vspace{0.2cm} 

(5.2) Let $X$ and $Y$ be smooth projective varieties. 
For an object ${\cal E} \in D(X \times Y)$, let 
$F_{{X \to Y};{{\cal E}}}: 
D(X) \to D(Y)$ be the Fourier-Mukai functor defined as 
$F_{{\cal E}}(\cdot) := {\bold R}(p_Y)_*({\bold L}(p_X)^*(\cdot)
\otimes^{\bold L}{\cal E})$. Here $p_X$ (resp. $p_Y$) is the 
first projection (resp. the second projection) of 
$X \times Y$. If there is an object ${\cal G} \in D(Y \times 
X)$ such that $F_{{Y \to X};{\cal G}}\circ 
F_{{X \to Y};{\cal E}} \cong id_X$ and 
$F_{{X \to Y};{\cal E}}\circ 
F_{{Y \to X};{\cal G}} \cong id_Y$, 
then $F_{{X \to Y};{\cal E}}$ 
is called an equivalence as a Fourier-Mukai transform. 
Note that the functor $\Psi$ in \S 4 is an equivalence 
as a Fourier-Mukai transform. 

Let $S$ be a (not necesarily projective) smooth algebraic 
variety. 

Let $p_{X \times Y \times S \to X 
\times Y}: X \times Y \times S \to 
X \times Y$ be the projection to the first and second 
factors (similarly, we define projections 
$p_{X \times Y \times S \to X \times S}$ and 
$p_{X \times Y \times S \to Y \times S}$).  
Put ${\cal E}\odot S := {\bold L}(p_{X \times Y \times S 
\to X \times Y})^*{\cal E} 
\in D(X \times Y \times S)$. 
Define a functor 

$F_{X \times S \to Y \times S;{\cal E}\odot S}: 
D(X \times S) \to D(Y \times S)$ 
as 
$${\bold R}
(p_{X \times Y \times S \to Y \times S})_*
({\bold L}(p_{X \times Y \times S \to X \times S})^*
(\cdot)\otimes ({\cal E}\odot S)).$$

Let $Z$ be a smooth projective variety.   
Let ${\cal G}$ be an object of $D(Y \times Z)$. 
By a similar argument 
as [Mu, Proposition 1.3], we have \vspace{0.12cm}

$F_{Y \times S \to 
Z \times S; {\cal G}\odot S} \circ 
F_{X \times S \to Y \times S;{\cal E}\odot S}$

$\cong F_{X \times S \to Z \times S; <{\cal G}, 
{\cal E}>
\odot S}.$  
\vspace{0.12cm}

Here 

$<{\cal G}, {\cal E}> := {\bold R}
(p_{X \times Y \times Z \to X \times Z})_*({\bold 
L}(p_{X \times Y \times Z \to X \times Y})^*{\cal E} 
\otimes  
{\bold L}(p_{X \times Y \times Z \to Y \times Z})^*
{\cal G}).$ 
\vspace{0.2cm}

{\bf Proposition (5.3)}. {\em Assume that $F_{X \to Y;{\cal E}}: 
D(X) \to D(Y)$ is an equivalence as a Fourier-Mukai transform. 
Then $F_{X \times S \to Y \times S; {\cal E}\odot S}: 
D(X \times S) \to D(Y \times S)$ is also an equivalence.}          
\vspace{0.15cm}

{\em Proof of (5.3)}. Since $F_{X \to Y;{\cal E}}$ 
is an equivalence as a Fourier-Mukai transform, one can 
find ${\cal G} \in D(Y \times X)$ such that  
$F_{{Y \to X};{\cal G}}\circ 
F_{{X \to Y};{\cal E}} \cong id_X$ and 
$F_{{X \to Y};{\cal E}}\circ 
F_{{Y \to X};{\cal G}} \cong id_Y$.  
This means that $F_{X \to X; <{\cal G}, {\cal E}>} 
\cong F_{X \to X, {\cal O}_{\Delta_X}}$ and 
$F_{Y \to Y; <{\cal E}, {\cal G}>} \cong 
F_{Y \to Y; {\cal O}_{\Delta_Y}}$. By a 
theorem of Orlov [O 2, Theorem 2.2], we conclude 
that $<{\cal G}, {\cal E}> \cong {\cal O}_{\Delta_X}$ 
in $D(X \times X)$ and 
$<{\cal E}, {\cal G}> \cong {\cal O}_{\Delta_Y}$ 
in $D(Y \times Y)$.  
Therefore, we have \vspace{0.12cm}

$ F_{Y \times S \to 
X \times S; {\cal G}\odot S} \circ 
F_{X \times S \to Y \times S;{\cal E}\odot S}$ 

$\cong F_{X \times S \to X \times S; <{\cal G}, 
{\cal E}>
\odot S}$ 

$\cong F_{X \times S \to X \times S; 
{\cal O}_{\Delta_X}
\odot S} \cong id_X.$   
\vspace{0.12cm}

Similarly, we have 

$$F_{X \times S \to 
Y \times S; {\cal E}\odot S} \circ 
F_{Y \times S \to X \times S;{\cal G}\odot S} 
\cong id_Y. $$  
\vspace{0.15cm}

(5.4)({\em Proof of} (5.1)): It is sufficient to 
show that, for $z, z' \in Z$ and for $i \in {\bold Z}$, 
${\mathrm Hom}^i({\cal O}_z, {\cal O}_{z'}) \cong 
{\mathrm Hom}^i(\Psi_Z({\cal O}_z), \Psi_Z({\cal O}_{z'}))$. 
Let $g: W \to S$ be the natural map induced by $f: Z \to 
{\bar Z}$.
We only have to consider the case where $g(z) 
= g(z')$. We put $s: = g(z)$. Let $s \in U(s)$ be an 
open set of $S$ such that $g^{-1}(U(s)) \cong 
{\bold P}^r \times U(s)$. Under this identification, 
let $(y,s)$ and $(y', s)$ be the points which 
correspond to $z$ and $z'$ respectively. 
Put $X:= {\bold P}
({\cal O}_{{\bold P}^r}\oplus \Theta_{{\bold P}^r})$. 
Then ${\bold P}^r \times U(s)$ can be identified with 
the (negative) section of 
$X \times U(s) \to {\bold P}^r \times U(s)$. Now, the 
Mukai flop $X - - \to X^+$ induces a birational map 
$X \times U(s) - - \to  X^+ \times U(s)$. Put 
$\hat X := X \times_{\bar X}X^+$ as in \S 1. Then we 
have an equivalence $F_{{\cal O}_{\hat X}\odot U(s)}: 
D(X \times U(s)) \to D(X^+ \times U(s))$ (cf. (5.3)). 
Note that, the formal completion of $Z$ along 
$g^{-1}(U(s))$ is isomorphic to the formal completion of 
$X \times U(s)$ along ${\bold P}^r \times U(s)$. By 
the same argument as (4.2) and (4.3), the homomorphism 
$${\mathrm Hom}^i_{D(Z)}({\cal O}_z, {\cal O}_{z'}) 
\to {\mathrm Hom}^i_{D(Z^+)}(\Psi_Z({\cal O}_z), 
\Psi_Z({\cal O}_{x'}))$$ 
can be identified with 
$${\mathrm Hom}^i_{D(X \times U(s))}({\cal O}_{(y,s)}, 
{\cal O}_{(y',s)}) 
\to {\mathrm Hom}^i_{D(X^+ \times U(s))}(F_{{\cal O}_{\hat X}
\odot U(s)}({\cal O}_{(y,s)}), F_{{\cal O}_{\hat X}
\odot U(s)}({\cal O}_{(y',s)})).$$ 
This homomorphism is an isomorphism 
by (5.3).  \vspace{0.15cm}

{\em Acknowledgement}:  The author thanks 
J. Sawon, B. Szendroi and R. Yoshizawa for useful comments. 
 
\vspace{0.2cm}

\begin{center}
Department of Mathematics, Graduate School of 
Sience, Osaka University, 

Toyonaka, Osaka 560, Japan 
\end{center}
 
\end{document}